\documentclass{amsart}

\theoremstyle{plain}
\newtheorem{lemma}{Lemma}[section]

\newtheorem{con}[lemma]{Conjecture}
\newtheorem{claim}{Claim}

\theoremstyle{definition}
\newtheorem{defi}[lemma]{Definition}

\newtheorem{example}[lemma]{Example}

\theoremstyle{remark}
\newtheorem{remark}[lemma]{Remark}
\newtheorem{notation}[lemma]{Notation}

\newcommand{\pic}{{\rm Pic}}
\newcommand{\bs}{{\rm Bs}}
\newcommand{\p}{\mathbb{P}}
\newcommand{\co}{\mathbb{C}}

\newcommand{\oc}{{\mathcal O}}
\newcommand{\lls}{\mathbb{L}}
\newcommand{\mms}{\mathbb{M}}

\newcommand{\ls}{{\mathcal L}}
\newcommand{\tls}{\tilde{\ls}}
\newcommand{\tl}{\tilde{L}}
\newcommand{\tf}{\tilde{F}}
\newcommand{\tm}{\tilde{M}}

\newcommand{\I}{{\mathcal I}}

\newcommand{\rt}{\rightarrow}

\numberwithin{equation}{section}
\numberwithin{figure}{section}
\numberwithin{table}{section}

\begin{document}

\title{A counterexample to a conjecture on linear systems on $\p^3$}
\author{Antonio Laface and Luca Ugaglia}
\address{Dipartimento di Matematica, Universit\`a degli Studi di Milano,
Via Saldini 50, 20133 MILANO}
\email{laface@mat.unimi.it\\
ugaglia@mat.unimi.it}
\begin{abstract}
In his paper \cite{ci} Ciliberto proposes a conjecture in order to
characterize special linear systems of $\p^n$ through multiple base points.
In this note we give a counterexample to this conjecture by showing that
there is a substantial difference between the speciality of linear systems
on $\p^2$ and those of $\p^3$.

\end{abstract}
\maketitle

\section*{Introduction}

Let us take the projective space $\p^n$ and
let us consider the linear system of hypersurfaces of degree $d$ having some points
of fixed multiplicity. The virtual dimension of such systems is
the dimension of the space of degree $d$ polynomials minus the conditions
imposed by the multiple points and the expected dimension is the maximum
between the virtual one and $-1$.
The systems whose dimension is bigger than the expected one are called
{\em special systems}.

There exists a conjecture due to Harbourne~\cite{har} and
Hirschowitz~\cite{hi}, characterizing special linear systems on
$\p^2$, which has been proved in some special cases
\cite{cm,cm2,mi,ev}.

Concerning linear systems on $\p^n$, in \cite{ci} Ciliberto gives a conjecture based
on the classification of special linear systems through double points.
In this note we describe a linear system on $\p^3$ that we found in a list of special systems
generated with the help of Singular and which turns out to be a
counterexample to that conjecture.

The paper is organized as follows: in Section 1 we fix some notations and state Ciliberto's conjecture, while Section 2 is devoted to the counterexample. In Section 3 we try to explain speciality of some systems by the Riemann-Roch formula, and we conclude the note with an appendix containing some computations.

\section{Preliminaries}

We start by fixing some notations.
\begin{notation}
Let us denote by $\lls_n(d,m_1^{a_1},\ldots,m_r^{a_r})$ the linear system of hypersurfaces of
$\p^n$ of degree $d$, passing through $a_i$ points with multiplicity $m_i$, for $i=1,\ldots,r$.
Let $\I_Z$ be the ideal of the zero dimensional scheme of multiple points. We denote by
$\ls_n(d,m_1^{a_1},\ldots,m_r^{a_r})$ the sheaf $\oc_{\p^n}(d)\otimes\I_Z$.\\
Given the system $\lls=\lls_n(d,m_1^{a_1},\ldots,m_r^{a_r})$, its {\em virtual dimension} is
\[
v(\lls)=\binom{d+n}{n}-\sum_{i=1}^ra_i\binom{m_i+n-1}{n}-1,
\]
and the {\em expected dimension} is
\[
e(\lls)=\max(v(\lls),-1).
\]

\end{notation}
A linear system will be called {\em special} if its expected dimension is strictly smaller
than the effective one.
\begin{remark}
Throughout the paper, if no confusion arises, we will use sometimes the same letter
to denote a linear system and the general divisor in the system.
\end{remark}
We recall the following definition, see \cite{ci}.

\begin{defi}
Let $X$ be a smooth, projective variety of dimension $n$, let $C$ be a
smooth, irreducible curve on $X$ and let $\mathcal{N}_{C|X}$ be the normal
bundle of $C$ in $X$. We will say that $C$ is a negative curve if there
is a line bundle $\mathcal{N}$ of negative degree and a surjective map
$\mathcal{N}_{C|X}\rt \mathcal{N}$.
The curve $C$ is called a $(-1)$-curve of size $a$, with $1\leq a\leq n-1$,
on $X$ if $C\cong \p^1$ and
$\mathcal{N}_{C|X}\cong \oc_{\p^1}(-1)^{\oplus a}\oplus\mathcal{N}$,
where $\mathcal{N}$ has no summands of negative degree.
\end{defi}

The main conjecture stated in \cite{ci} is the following.

\begin{con}
\label{cil}
Let $X$ be the blow-up of $\p^n$ at general points $p_1,\ldots,p_r$ and let
$\lls = \lls_n(d,m_1,\ldots,m_r)$ be a linear system with multiple base points
at $p_1,\ldots,p_r$.
Then:
\begin{itemize}
\item[(i)] the only negative curves on $X$ are $(-1)$-curves;
\item[(ii)] $\lls$ is special if and only if there is a $(-1)$-curve $C$ on $X$
corresponding to a curve $\Gamma$ on $\p^n$ containing $p_1\ldots,p_r$ such
that the general member $D\in\lls$ is singular along $\Gamma$;
\item[(iii)] if $\lls$ is special, let $B$ be the component of the base locus
of $\lls$ containing $\Gamma$ according to Bertini's theorem. Then the
codimension of $B$ in $\p^n$ is equal to the size of $C$ and $B$ appears
multiply in the base locus scheme of $\lls$.
\end{itemize}
\end{con}

In this note we give a counterexample to points (ii) and (iii) of this
conjecture.

\section{Counterexample}

Let us consider the linear system of surfaces of degree nine with one point of
degree six and eight points of degree four in $\p^3$, i.e. the system
$\lls = \lls_3(9,6,4^8)$. In this section we are going to study this system, showing
in particular that it is special but its general member is not singular along a rational curve.\\
If we denote by $Q = \lls_3(2,1,1^8)$ the quadric through the nine simple points,
we have the following:
\begin{claim}
$\lls_3(9,6,4^8) = Q + \lls_3(7,5,3^8).$
\end{claim}
If we denote by $H_1,H_2$ two generators of $\pic(Q)$,
considering the restriction $\lls_{|Q}$ we get the system of curves in $|9H_1+9H_2|$, with
one point of multiplicity $6$ and eight points of multiplicity $4$.
We denote for short this system by  $|9H_1+9H_2|-6p_0-\sum 4p_i$.\\
Looking the Appendix \ref{app-1}, we can see that $|9H_1+9H_2|-6p_0-\sum 4p_i$ corresponds to
the planar system $\lls_2(12,3^2,4^8)$. This last system can not be
$(-1)$-special (see the Appendix \ref{app-2}) and $v(\lls_2(12,3^2,4^8)) = -2$. Therefore, by \cite{mi}
we may conclude that it is empty. \\
In particular, also $\lls_{|Q} = \emptyset$, and hence $\lls$ must contain
$Q$ as a fixed component. By subtracting $Q$ from $\lls$ we get $\lls_3(7,5,3^8)$, which
proves our claim.
\vskip .2truecm

This means that the free part of $\lls$ is contained in $\lls_3(7,5,3^8)$ which has
virtual dimension $4$. So $\lls$ is a special system. \\
In order to show that $\lls$ gives a counterexample to Conjecture \ref{cil} we are now
going to prove that the general member of $\lls$ is singular only along the curve $C$,
intersection of $Q$ and $\lls_3(7,5,3^8)$, and that $C$ does not contain
rational components.\\
We can consider $C$ as the restriction $\lls_3(7,5,3^8)_{|Q}$.
This is equal to $|7H_1+7H_2|-5p_0-\sum 3p_i$ on the quadric $Q$, which
corresponds to $\lls_2(9,2^2,3^8)$ on $\p^2$. This system
is not special of dimension $0$ and it does not contain rational components
(see Appendix \ref{app-3}). \\

Clearly the curve $C$ is contained in $\lls_{\rm sing}$ (i.e. the singular locus
of $\lls$). We are going to show that in fact $C = \lls_{\rm sing}$. \\
First of all, let us denote by $\lls_3(7,5,3^8,1_Q)$ the subsystem
of $\lls_3(7,5,3^8)$ obtained by imposing one general simple point on the quadric.
Since $\lls_3(7,5,3^8,1_Q)_{|Q}=\emptyset$, $Q$ is a fixed component of this system
and the residual part is given by $\lls_3(5,4,2^8)$. Now $\lls_3(5,4,2^8)_{|Q}$
is the system $|5H_1+5H_2|-4p_0-\sum 2p_i$ which corresponds to the non-special system
$\lls_2(6,1^2,2^8)$, of dimension $1$. Therefore, imposing two general simple points on $Q$
and restricting we get that the system $\lls_3(5,4,2^8,1^2_Q)_{|Q}$ is empty,
which implies that $\lls_3(5,4,2^8,1^2_Q)$ has $Q$ as a fixed component.
The residual system $\lls_3(3,3,1^8)$
is non-special of dimension $1$ (because each surface of this system is a cone over
a plane cubic through eight fixed points). This implies that the effective dimension of
$\lls_3(5,4,2^8)$ can not be greater than $3$. Therefore it must be $3$ since the virtual
dimension is $3$. By the same argument one shows that the effective dimension of
$\lls_3(7,5,3^8)$ is $4$. \\
Observe that $\bs(\lls_3(7,5,3^8))\subseteq\bs(2Q + \lls_3(3,3,1^8))$ since
$2Q + \lls_3(3,3,1^8)\subseteq\lls_3(7,5,3^8)$. So $\bs(\lls_3(7,5,3^8))$ could have only $Q$
as fixed component, but this is not the case since
$\dim \lls_3(7,5,3^8) = \dim \lls_3(5,4,2^8) + 1$. The only curves that may belong to
$\bs(\lls_3(7,5,3^8))$ are the genus $2$ curve $C=\lls_3(7,5,3^8)_{|Q}$
and the nine lines of $\bs(\lls_3(3,3,1^8))$ trough the vertex of the cone and each one of the
nine base points of the pencil of plane cubics. \\
We can then conclude that the singular locus $\lls_{\rm sing}$
consists only of the curve $C$, since the subsystem $3Q + \lls_3(3,3,1^8)$ is not
singular along the nine fixed lines.

\section{Speciality and Riemann-Roch theorem}
Let $Z$ be a zero dimensional scheme of $\p^3$ and $\I_Z$ be its ideal sheaf.
We put $\ls=\oc_{\p^3}(d)\otimes\I_Z$ and consider the exact sequence
\begin{eqnarray*}
0\rt H^0(\ls)\rt H^0(\oc_{\p^3}(d))\rt H^0(\oc_Z)
\rt H^1(\ls)\rt H^1(\oc_{\p^3}(d))\rt H^1(\oc_Z) \\
\rt H^2(\ls)\rt H^2(\oc_{\p^3}(d))\rt H^2(\oc_Z)
\rt H^3(\ls)\rt H^3(\oc_{\p^3}(d))\rt H^3(\oc_Z),
\end{eqnarray*}
obtained tensoring by $\oc_{\p^3}(d)$ the sequence defining $Z$ and taking cohomology.
From this sequence we obtain that $h^i(\ls) = h^i(\oc_{\p^3}(d)) = 0$ for
$i=2,3$ since $h^i(\oc_Z)=0$ for $i=1,2,3$. We also obtain that the virtual dimension of $\lls$,
$h^0(\oc_{\p^3}(d))-h^0(\oc_z)-1$ is equal to $h^0(\ls)-h^1(\ls)-1$ and hence to $\chi(\ls)-1$. \\
If $Z=\sum m_ip_i$ is a scheme of fat points, then on the blow-up $X\stackrel{\pi}{\rt}\p^3$
along these points we may consider the divisor $\tl=\pi^*\oc_{\p^3}(d)-\sum m_iE_i$ and the associated sheaf $\tls=\oc_X(\tl)$.
Since $h^i(X,\tls)=h^i(\p^3,\ls)$, the virtual dimension of $\lls$ is equal to $\chi(\tls)-1$.
By Riemann-Roch formula (see \cite{ha}) for a divisor $\tl$ on the threefold $X$,
\[
\chi(\tl) = \frac{\tl(\tl-K_X)(2\tl-K_X)+c_2(X)\cdot\tl}{12} + \chi(\oc_X),
\]
we obtain the following formula for the virtual dimension of $\lls$:
\[
v(\lls) = \frac{\tl(\tl-K_X)(2\tl-K_X)+c_2(X)\cdot\tl}{12}
\]
since $\chi(\oc_X)=1$. \\
If the linear system $\lls$ can be written as $F+\mms$, where $F$ is the fixed divisor and
$\mms$ is a free part, then on $X$ we have $|\tl| = \tf + |\tm|$.
Therefore the above formula says that
\[
v(\tl) = v(\tf) + v(\tm) + \frac{\tf\tm(\tl-K_X)}{2}.
\]
Let us suppose that the residual system $\mms$ is non-special.
The system $\lls$ has the same effective dimension as $\mms$, while their virtual dimensions
differ by $v(\tf) + \tf\tm(\tl-K_X)/2$. Therefore we can conclude that
$\lls$ is special if $v(\tf) + \tf\tm(\tl-K_X)/2$ is smaller than zero.
\begin{example}
For instance, let us consider the system $\lls:=\lls_3(4,2^9)$. It is special because
its virtual dimension is $-2$ while it is not empty since it is equal to $2Q$, where $Q$
is the quadric through the nine simple points. In this case $F=2Q$ and $\mms=\co$, so
$v(F)=-2$ and $\tf\tm(\tl-K_X)/2=0$.
\end{example}
\begin{example}
Let us consider now the example we described in the previous section, i.e. the
system $\lls_3(9,6,4^8)$. We have seen that it can be written as $Q+\mms$,
where $Q$ is the quadric through the nine points, while $\mms=\lls_3(7,5,3^8)$
is the residual free part. The Chow ring $A^*(X)$ (where $X$ is the blow-up of $\p^3$
along the nine simple points) is generated by $\langle H,E_0,E_1,\ldots,E_8\rangle$, where
$H$ is the pull-back of the hyperplane divisor of $\p^3$ and the $E_i$'s are the exceptional
divisors. The second Chow group $A^2(X)$ is generated by $\langle h,e_0,e_1,\ldots,e_8\rangle$,
where $h=H^2$ is the pull-back of a line, while $e_i=-E_i^2$ is the class of a line inside
$E_i$, for $i=0,1,\ldots,8$. Clearly $H\cdot E_i=E_i\cdot E_j=0$ for $i\neq j$.
With this notation we can write:
\begin{eqnarray*}
\lls & = & |9H-6E_0-\sum 4E_i| \\
\mms & = & |7H-5E_0-\sum 3E_i| \\
Q   & = & 2H-E_0-\sum E_i \\
K_X & = & -4H+2E_0+\sum 2E_i.
\end{eqnarray*}
Therefore $Q\cdot\mms=14h-5e_0-\sum 3e_i$,  $\ls-K_X=13H-8E_0-\sum 6E_i$ and hence
$Q\tm(\tl-K_X)/2=-1$ (while $v(Q)=0$), which implies the speciality of $\lls$.
\end{example}

\section{Appendix}

\subsection{Linear systems on a quadric}
\label{app-1}\ \\

In order to study linear systems on a quadric $Q$ it may be helpful to transform them into planar systems by mean of a birational transformation $Q\rt\p^2$ obtained by blowing up a point and contracting the strict transforms of the two lines through it.
Such transformation gives rise to a $1:1$ correspondence between linear systems with one multiple point on the quadric and linear systems with two multiple points on $\p^2$.\\
In fact, let us consider a linear system $|aH_1+bH_2|-mp$ (i.e. a system of curves of kind
$(a,b)$ through one point $p$ of multiplicity $m$). Blowing up at $p$, one obtains
the complete system $|a\pi^*H_1+b\pi^*H_2-mE|$ which may be written as
$|(a+b-m)(\pi^*H_1+\pi^*H_2-E)-(b-m)(\pi^*H_1-E)-(a-m)(\pi^*H_2-E)|$.
Since the divisors $\pi^*H_i-E$ ($i=1,2$) are $(-1)$-curves, they may be contracted
giving a linear system on $\p^2$ of degree $a+b-m$ through two points of multiplicity
$b-m$ and $a-m$ and hence
\[
|aH_1+bH_2|-mp \rt \lls_2(a+b-m,b-m,a-m).
\]

\vskip .5truecm

\subsection{$(-1)$-curves}
\label{app-2}\ \\

In order to study the speciality of the systems
$\lls_2(12,3^2,4^8),\  \lls_2(9,2^2,3^8)$ and $\lls_2(6,1^2,2^8)$,
we need to produce a complete list
of all the $(-1)$-curves of $\p^2$ of kind $\lls_2(d,m_1,m_2,m_3,\ldots,m_{10})$ which
may have an intersection less than -1 with some of these systems.
Clearly it is enough to consider the system $\lls_2(12,3^2,4^8)$, whose degree and multiplicities are the biggest.
From the condition of being contained twice in this system we deduce the following inequalities:
$d\leq 6,\  0\leq m_1,m_2\leq 1$ and $0\leq m_3,\ldots,m_{10}\leq 2$. Moreover let us see that
$m_3=\ldots=m_{10}=m$. Otherwise the system would contain twice the compound $(-1)$-curve given
by the union of all the simple $(-1)$-curves obtained by permuting the points
$p_3,\ldots,p_{10}$. In this case the multiplicities of the compound curve at these points
would be too big. An explicit calculation shows that the only $(-1)$-curve of the form
$\lls_2(d,m_1,m_2,m^8)$ satisfying the preceding conditions is $\lls_2(1,1,1,0^8)$, but
this has non negative intersection with any of these systems.

\vskip .5truecm

\subsection{$\lls_2(9,2^2,3^8)$ does not contain rational components}
\label{app-3}\ \\

Let $S$ be the blow up of $\p^2$ along the ten points and let $C$ be the strict transform
of the curve given by $\lls=\lls_2(9,2^2,3^8)$.
Suppose that there exists an irreducible rational component $C_1$ of $C$. Observe that $v(C_1)=0$ since the system $|C_1|$ has dimension $0$ and it is non-special by \cite{mi}. Therefore, from
$g(C_1)=v(C_1)=0$, we get that $C_1$ is a $(-1)$-curve.

We are going to see that if this is the case, then $C\cdot C_1=-1$. Let us take the
following exact sequence:
\[
0\rt \oc_S(C-C_1) \rt \oc_S(C) \rt \oc_{\p^1}(C\cdot C_1)\rt 0.
\]
By the subsection above, $h^1(\oc_S(C))=0$. Let us see that also $h^1(\oc_S(C-C_1))=0$. Otherwise the system $|C-C_1|$ would be special and in particular, by \cite{mi} there would exist a $(-1)$-curve $C_2$ such that $C_2\cdot (C-C_1)\leq -2$. Since $\lls$ is non-special,
$C\cdot C_2\geq -1$ and hence $C_1\cdot C_2\geq 1$. This implies that $|C_1+C_2|$ has
dimension at least $1$, which is impossible since $C_1+C_2$ is contained in the fixed
locus of $\lls$. Since $h^0(\oc_S(C-C_1))=h^0(\oc_S(C))=1$, the cohomology of the preceding sequence gives $h^0(\oc_{\p^1}(C\cdot C_1))=h^1(\oc_{\p^1}(C\cdot C_1))=0$, which means that $C\cdot C_1 = -1$ as claimed before.

Arguing as in the previous subsection, we get $|C_1|=\lls_2(d,\- m_1,\- m_2\- ,m^8)$
with $d\leq 9,\ 0\leq m_1,m_2\leq 2,\ 0\leq m\leq 3$. An easy computation shows that the only $(-1)$-curve of this form is $\lls_2(1,1,1,0^8)$ and in this case $C\cdot C_1 = 5$.

\end{document}